\documentclass[12pt]{article}

\usepackage{amssymb}
\usepackage{amsmath,amsthm}
\usepackage[dvips]{graphicx}
\usepackage{wrapfig}
\usepackage{url}

\newtheorem{theorem}{Theorem}[section]

\newtheorem{lemma}[theorem]{Lemma}
\newtheorem{proposition}[theorem]{Proposition}
\newtheorem{question}[theorem]{Question}
\theoremstyle{definition}

\theoremstyle{remark}

\title{The reductivity of spherical curves}
\author{ Ayaka Shimizu \thanks{Department of Mathematics, Gunma National College of Technology, 580 Toriba-cho, Maebashi-shi, Gunma, 371-8530, Japan.  shimizu1984@gmail.com}}
\date{\today}

\begin{document}

\maketitle

\begin{abstract}
We show that we can obtain a reducible spherical curve from any non-trivial spherical curve by four or less inverse-half-twisted splices, i.e., the reductivity, which represents how reduced a spherical curve is, is four or less. 
We also discuss unavoidable sets of tangles for spherical curves. 
\end{abstract}

\section{Introduction}

\begin{figure}[ht]
 \begin{center}
  \includegraphics[width=100mm]{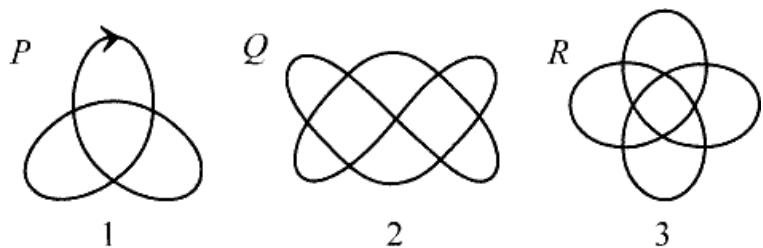}
 \end{center}
 \caption{Spherical curves.}
 \label{ex-reductivity}
\end{figure}

\begin{figure}[ht]
 \begin{center}
  \includegraphics[width=40mm]{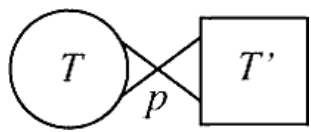}
 \end{center}
 \caption{Reducible spherical curve.}
 \label{red}
\end{figure}

A {\it spherical curve} is a smooth immersion of the circle into the sphere where the self-intersections are transverse and double points (we call the double point {\it crossing}). 
In this paper we assume every spherical curve is oriented, and has at least one crossing. 
We represent, if necessary, the orientation of a spherical curve by an arrow as depicted in the left-hand side of Fig. \ref{ex-reductivity}. 
A spherical curve $P$ is \textit{reducible} and has a \textit{reducible crossing} $p$ if $P$ has a crossing $p$ as shown in Fig. \ref{red}, where $T$ and $T'$ imply parts of the spherical curve. 
$P$ is \textit{reduced} if $P$ is not reducible such as the spherical curves in Fig. \ref{ex-reductivity}. 
Note that around a reducible (resp. non-reducible) crossing, there are exactly three (resp. four) disjoint regions, where a {\it region} of a spherical curve is a part of the sphere divided by the spherical curve. \\

\begin{figure}[ht]
 \begin{center}
  \includegraphics[width=100mm]{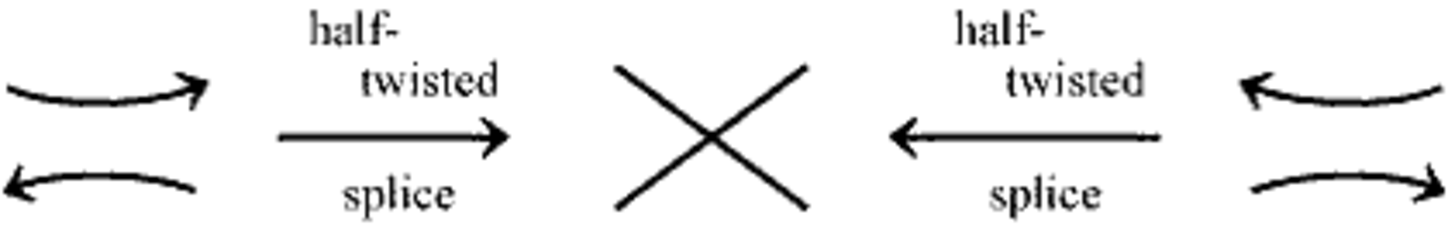}
 \end{center}
 \caption{Half-twisted splice.}
 \label{halftwisted}
\end{figure}

\begin{figure}[ht]
 \begin{center}
  \includegraphics[width=100mm]{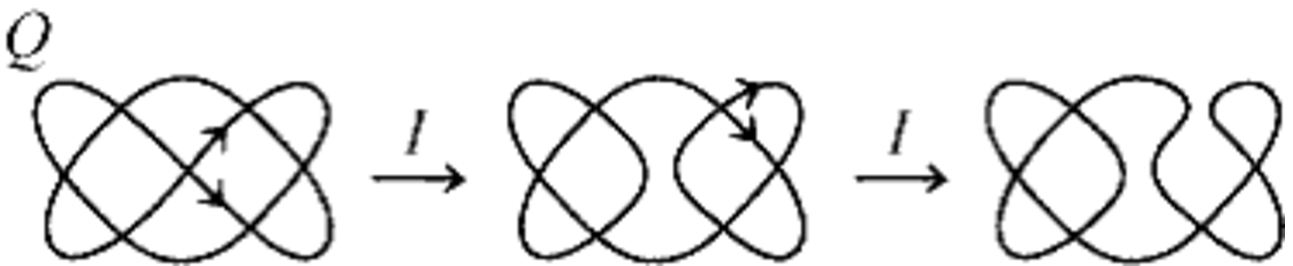}
 \end{center}
 \caption{The reductivity of $Q$ is 2.}
 \label{eight-two}
\end{figure}

A {\it half-twisted splice} is the local transformation on spherical curves as depicted in Fig. \ref{halftwisted} (\cite{calvo, IS}). 
\begin{figure}[ht]
 \begin{center}
  \includegraphics[width=55mm]{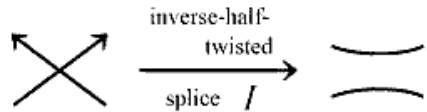}
 \end{center}
 \caption{Inverse-half-twisted splice.}
 \label{halftwisted-i}
\end{figure}
Then the inverse is the transformation depicted in Fig. \ref{halftwisted-i}. 
In this paper we call the inverse of the half-twisted splice {\it inverse-half-twisted splice}, and denote by $I$. 
We remark that the half-twisted splice and the inverse-half-twisted splice do not preserve the orientation of spherical curves. 
Then we give an orientation again to the spherical curve we obtain. 
We also remark that the half-twisted splice and the inverse-half-twisted splice do not depend on the orientations of spherical curves, but depend only on the relative orientations. 
Now we define the reductivity. 
The {\it reductivity} $r(P)$ of a spherical curve $P$ is the minimal number of $I$ which are needed to obtain a reducible spherical curve from $P$. 
For example, a reducible spherical curve has the reductivity 0, and the spherical curves $P$, $Q$ and $R$ in Fig. \ref{ex-reductivity} have the reductivity 1, 2 and 3, respectively (see Fig. \ref{eight-two} for $Q$. 
Note that we can not obtain a reducible curve by any single $I$ from $Q$.) 
In this paper, we show the following theorem:

\phantom{x}
\begin{theorem}
Every spherical curve has the reductivity four or less. 
\label{main}
\end{theorem}
\phantom{x}

\noindent This implies that we can obtain a reducible spherical curve from any spherical curve by four or less $I$. 
We have the following question. 

\phantom{x}
\begin{question}
Is it true that every spherical curve has the reductivity three or less?
\label{main-q}
\end{question}
\phantom{x}

\noindent In other words, is it true that there are no spherical curve with reductivity four?
The rest of this paper is organized as follows: 
In Section 2, we discuss the properties of reductivity by considering chord diagrams, and prove Theorem \ref{main}. 
In Section 3, we study the unavoidable sets of tangles for spherical curves as an approach to Question \ref{main-q}.

\section{Proof of Theorem\ref{main}}

In this section we show Theorem \ref{main} by using chord diagrams. 
We consider a spherical curve $P$ as an immersion $P: S^1 \rightarrow S^2$ of the circle into the sphere with some double points (crossings). 
A {\it chord diagram} for a spherical curve is an oriented circle considered as the preimage of the immersed circle with chords connecting the preimages of each crossing (\cite{GPV}).
\begin{figure}[ht]
 \begin{center}
  \includegraphics[width=60mm]{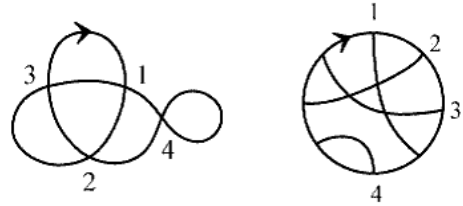}
 \end{center}
 \caption{Chord diagram.}
 \label{chord-ex}
\end{figure}

\begin{figure}[ht]
 \begin{center}
  \includegraphics[width=130mm]{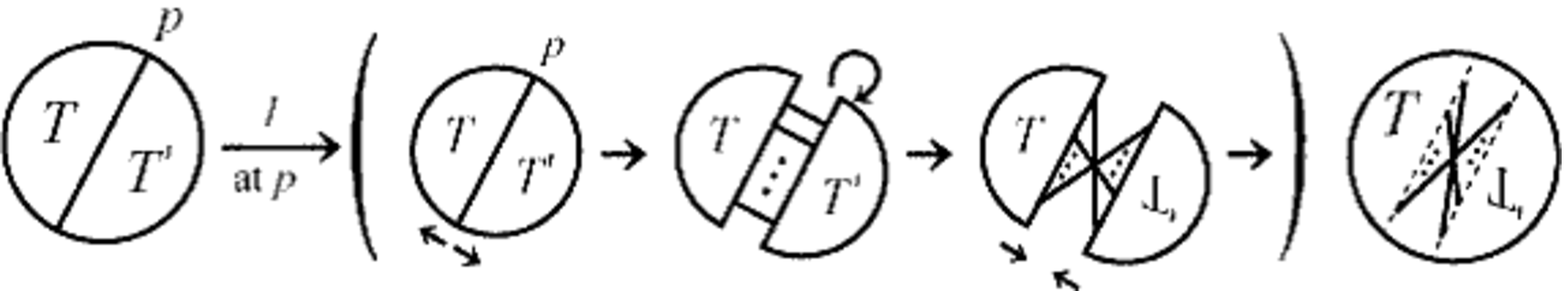}
 \end{center}
 \caption{$I$ on a chord diagram.}
 \label{a-on-chord}
\end{figure}

\begin{figure}[ht]
 \begin{center}
  \includegraphics[width=80mm]{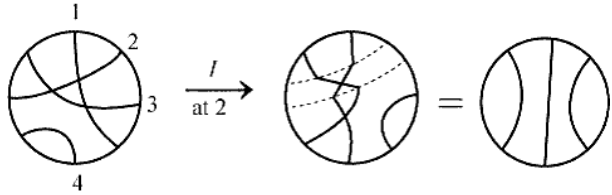}
 \end{center}
 \caption{$I$ at 2 on a chord diagram.}
 \label{a-on-chord-ex}
\end{figure}

A chord diagram is for a reducible spherical curve if and only if the chord diagram could have a chord without crossings. 
For example, the chord diagram in Fig. \ref{chord-ex} has the chord labeled 4 without crossings. 
The move $I$ at a crossing $p$ on a spherical curve corresponds to the move shown in Fig. \ref{a-on-chord} on the chord diagram; 
cut the chord diagram along the chord labeled $p$ (then $p$ vanishes), turn over the semicircle, and join the two semicircles again. 
For example, $I$ at the crossing labeled $2$ on the spherical curve in Fig. \ref{chord-ex} is represented on the chord diagram as shown in Fig. \ref{a-on-chord-ex}. 

\begin{figure}[ht]
 \begin{center}
  \includegraphics[width=70mm]{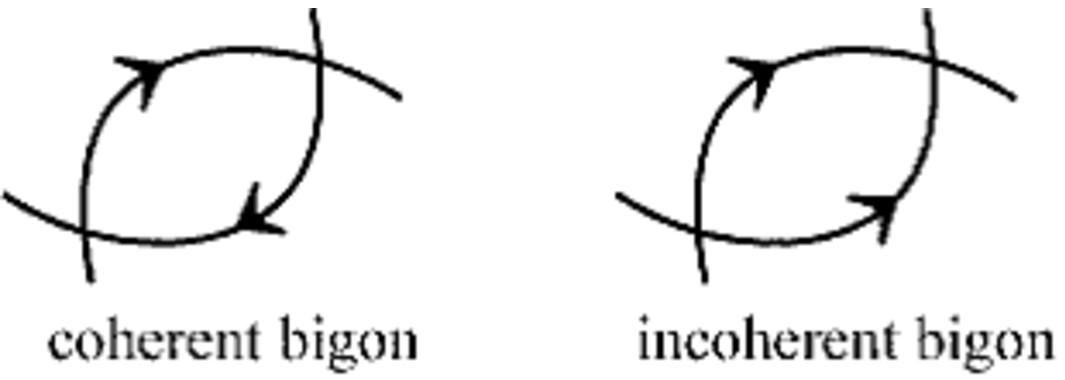}
 \end{center}
 \caption{Coherent and incoherent bigons.}
 \label{bigons}
\end{figure}

\begin{figure}[ht]
 \begin{center}
  \includegraphics[width=75mm]{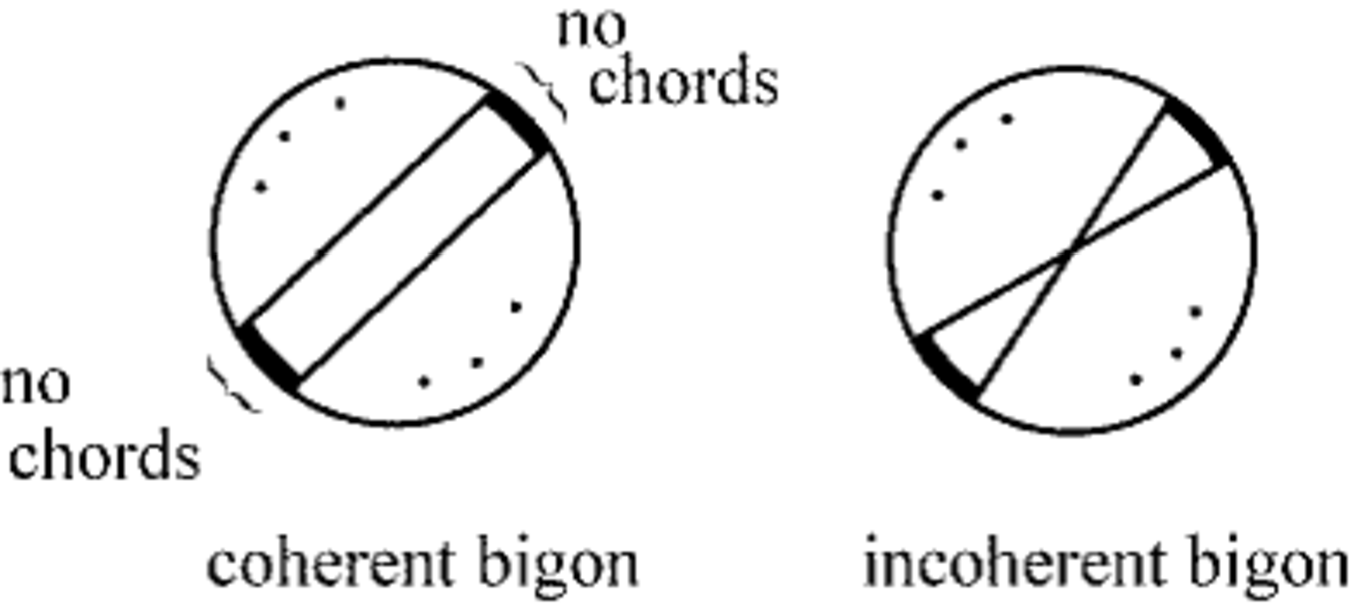}
 \end{center}
 \caption{Coherent and incoherent bigons in chord diagrams.}
 \label{chord-bigon}
\end{figure}

A region of a spherical curve is {\it coherent} (resp. {\it incoherent}) if the boundary of the region is coherent (resp. incoherent) (see, for example, Fig. \ref{bigons}). 
A coherent bigon and an incoherent bigon are represented by chord diagrams as shown in Fig. \ref{chord-bigon}. 
For coherent and incoherent bigons, we have the following lemmas:

\phantom{x}
\begin{lemma}
If a spherical curve $P$ has an incoherent bigon, then $r(P)\le 1$. 
\begin{proof}
By applying $I$ at one of the crossing on an incoherent bigon, we obtain a reducible spherical curve as shown in Fig. \ref{c-non-c-2}. 
\begin{figure}[ht]
 \begin{center}
  \includegraphics[width=60mm]{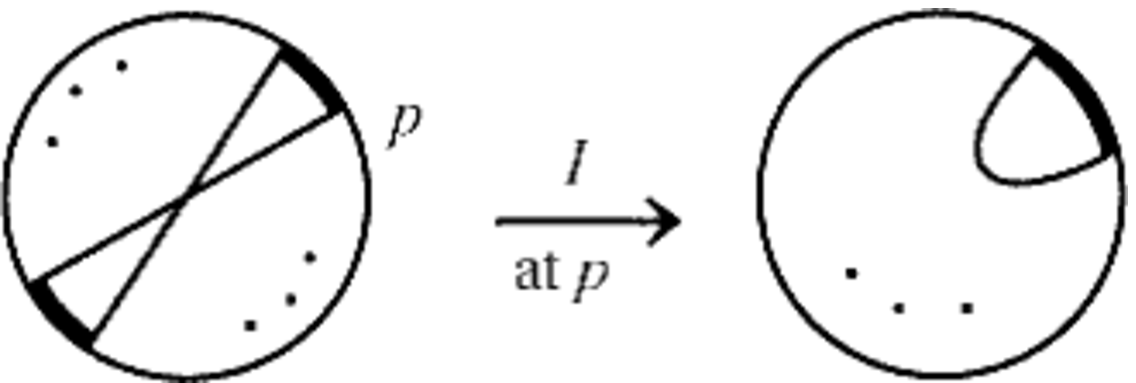}
 \end{center}
 \caption{Incoherent bigon.}
 \label{c-non-c-2}
\end{figure}
\end{proof}
\label{non-c-2}
\end{lemma}

\phantom{x}
\begin{lemma}
If a spherical curve $P$ has a coherent bigon, then $r(P)\le 2$. 
\begin{proof}
If $P$ is reducible, the reductivity is zero. 
If $P$ is reduced, there is a crossing $p$ as shown in the chord diagram in Fig. \ref{c-c-2}. 
By applying $I$ at $p$, we obtain a spherical curve which has an incoherent bigon. 
\begin{figure}[ht]
 \begin{center}
  \includegraphics[width=60mm]{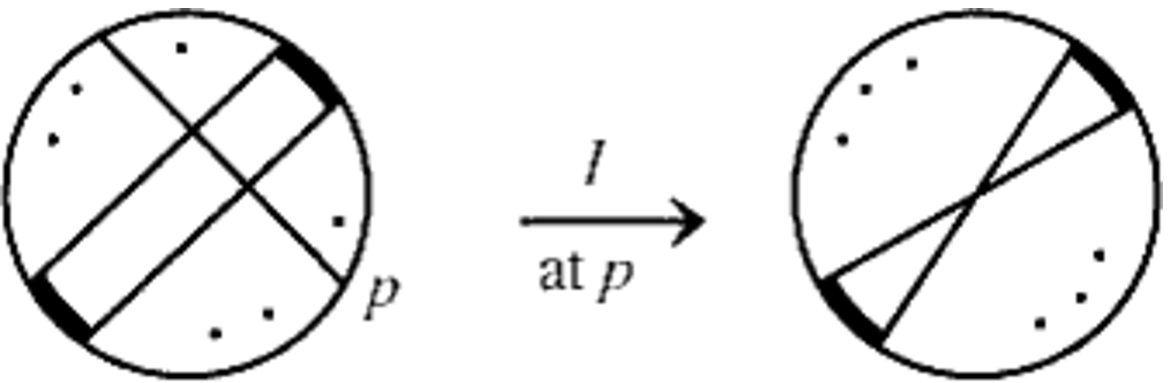}
 \end{center}
 \caption{Coherent bigon.}
 \label{c-c-2}
\end{figure}
Hence $P$ has the reductivity two or less. 
\end{proof}
\label{c-2}
\end{lemma}
\phantom{x}

\begin{figure}[ht]
 \begin{center}
  \includegraphics[width=100mm]{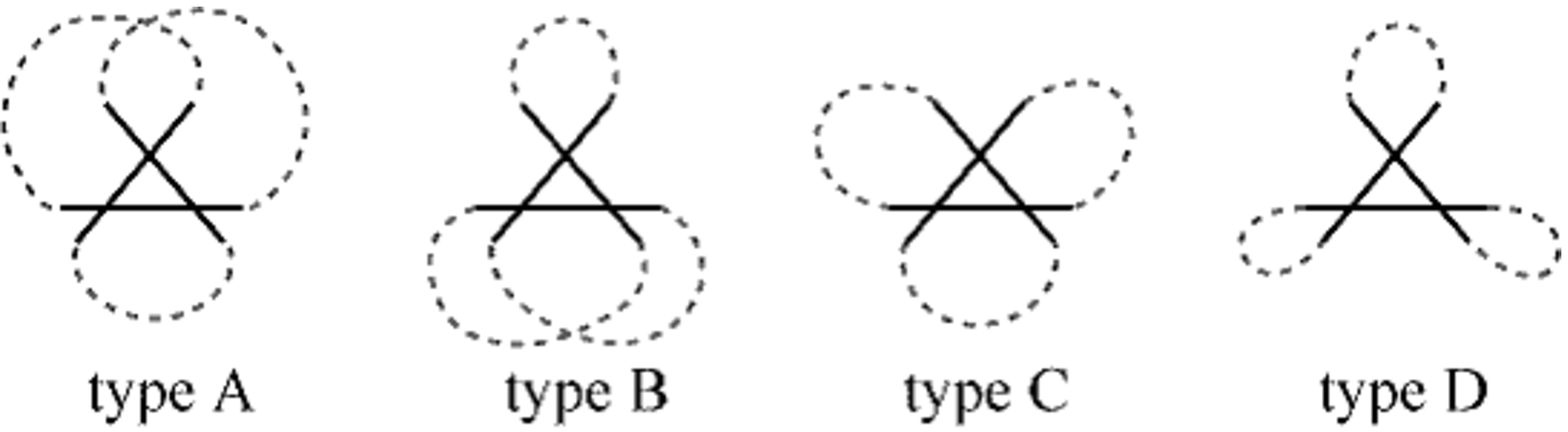}
 \end{center}
 \caption{The four types of trigons.}
 \label{abcd}
\end{figure}

\begin{figure}[ht]
 \begin{center}
  \includegraphics[width=100mm]{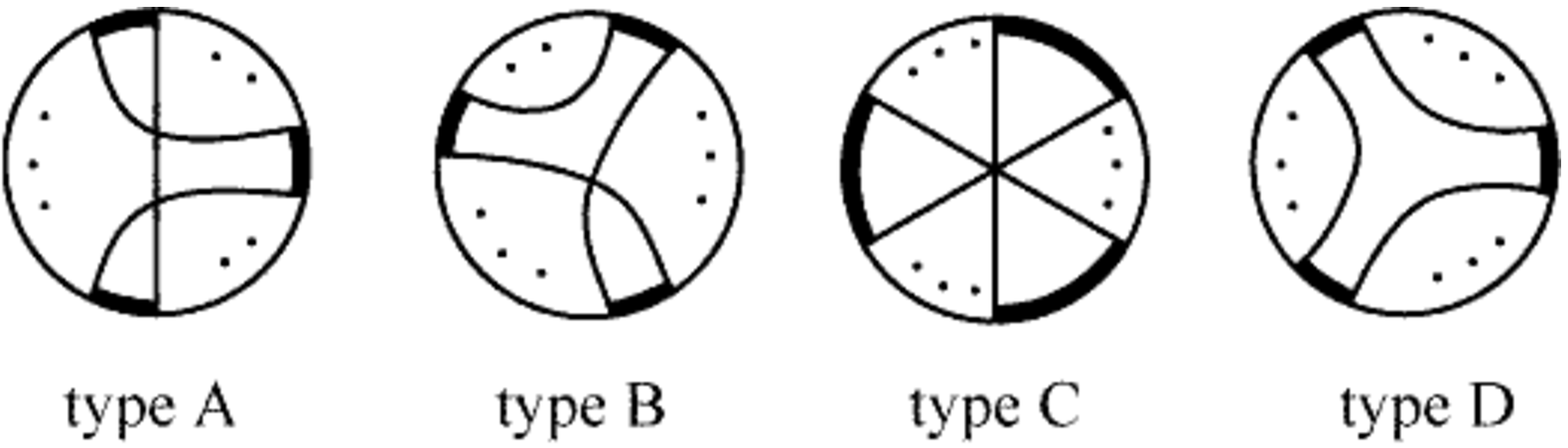}
 \end{center}
 \caption{The four types of trigons in chord diagrams.}
 \label{chord-trigon}
\end{figure}

A trigon of a spherical curve is one of the types A, B, C and D in Fig. \ref{abcd} with respect to the outer connections. 
These trigons are represented in chord diagrams as shown in Fig. \ref{chord-trigon}. 
We have the following lemmas for trigons of type A and B.

\phantom{x}
\begin{lemma}
If a spherical curve $P$ has a trigon of type A, then $r(P)\le 2$. 
\begin{proof}
By applying $I$ at $p$ in Fig. \ref{c-3a}, we have a spherical curve which has an incoherent bigon. 
\begin{figure}[ht]
 \begin{center}
  \includegraphics[width=60mm]{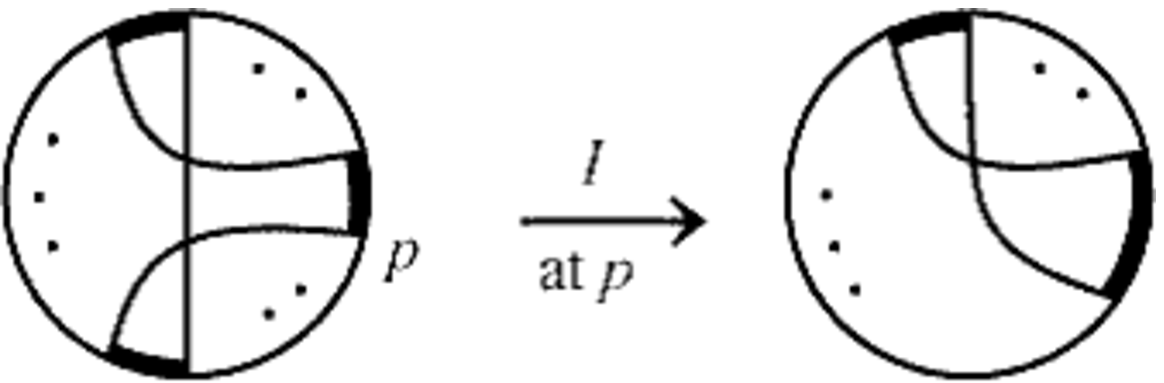}
 \end{center}
 \caption{Trigon of type A.}
 \label{c-3a}
\end{figure}
\end{proof}
\label{3a}
\end{lemma}
\phantom{x}

\begin{lemma}
If a spherical curve $P$ has a trigon of type B, then $r(P)\le 3$. 
\begin{proof}
By applying $I$ at $p$ in Fig. \ref{c-3b}, we have a spherical curve which has a coherent bigon. 
\begin{figure}[ht]
 \begin{center}
  \includegraphics[width=60mm]{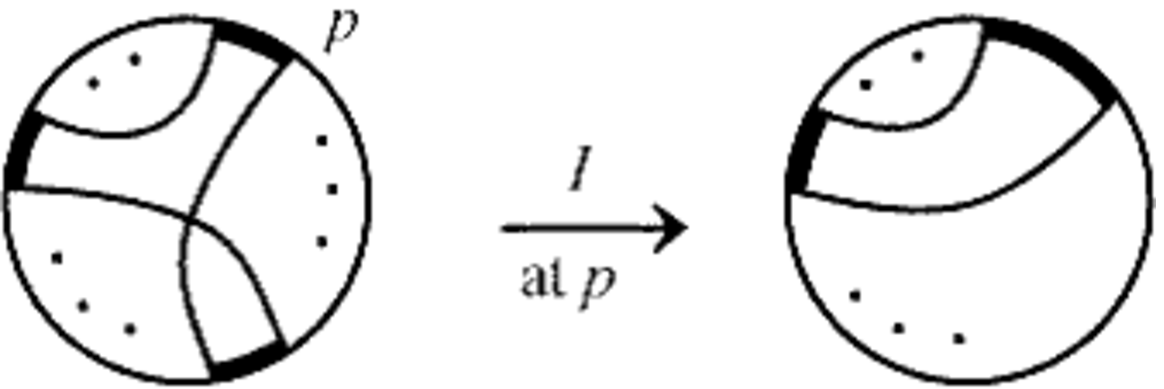}
 \end{center}
 \caption{Trigon of type B.}
 \label{c-3b}
\end{figure}
\end{proof}
\label{3b}
\end{lemma}
\phantom{x}

\begin{figure}[ht]
 \begin{center}
  \includegraphics[width=90mm]{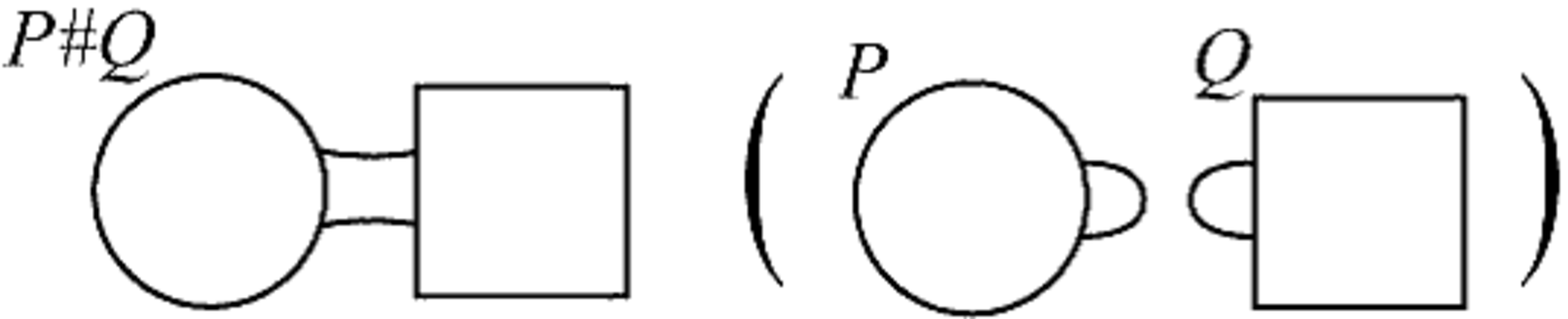}
 \end{center}
 \caption{Connected sum.}
 \label{conn}
\end{figure}

\noindent A \textit{connected sum} $P\sharp Q$ of two spherical curves $P$ and $Q$ is a spherical curve as depicted in Fig. \ref{conn}. 
For trigons of type C,  we have the following lemma:

\phantom{x}
\begin{lemma}
If a spherical curve $P$ has a trigon of type C, then $r(P)\le 3$. 
\begin{proof}
If the arcs $l$, $m$, and $n$ in Fig. \ref{lmn-c} around a trigon of type C have no mutual crossings, then $P$ is a connected sum of a trefoil ({\it trefoil} is the spherical curve depicted in the left-hand side of Fig. \ref{ex-reductivity}) and some (trivial or non-trivial) spherical curves as depicted in Fig. \ref{conn-trefoil}, whose reductivity is one or zero. 
\begin{figure}[ht]
 \begin{center}
  \includegraphics[width=30mm]{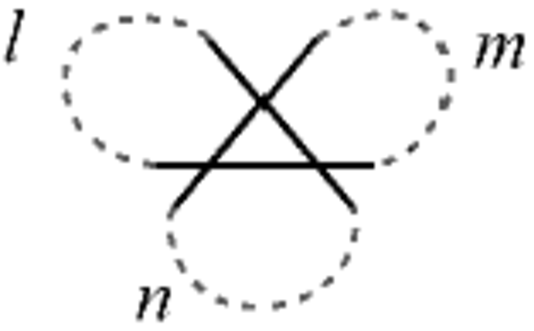}
 \end{center}
 \caption{Arcs of a trigon of type C.}
 \label{lmn-c}
\end{figure}
\begin{figure}[ht]
 \begin{center}
  \includegraphics[width=30mm]{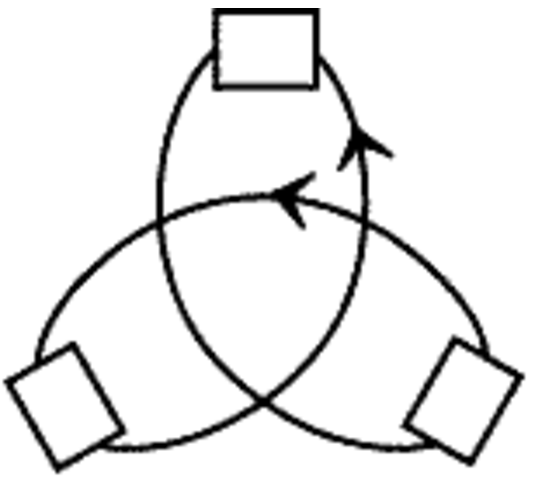}
 \end{center}
 \caption{Connected sums to a trefoil.}
 \label{conn-trefoil}
\end{figure}
If $P$ has a mutual crossing $p$ of $l$, $m$ and $n$, we obtain a spherical curve which has a trigon of type A by applying $I$ at $p$ (see Fig. \ref{c-3c}). 
\begin{figure}[ht]
 \begin{center}
  \includegraphics[width=60mm]{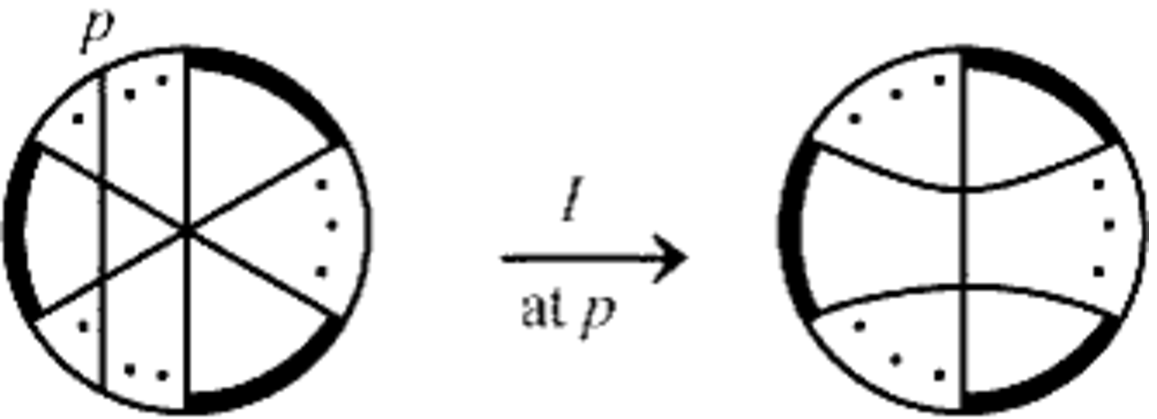}
 \end{center}
 \caption{Trigon of type C.}
 \label{c-3c}
\end{figure}
\end{proof}
\label{3c}
\end{lemma}
\phantom{x}

\noindent For type D, we have the following lemma:

\phantom{x}
\begin{lemma}
If a spherical curve $P$ has a trigon of type D, then $r(P)\le 4$. 
\begin{proof}
If $P$ is reducible, the reductivity is zero. 
If $P$ is reduced, $P$ has a crossing $p$ as shown in the chord diagram in Fig. \ref{c-3d}. 
Apply $I$ at $p$, and we obtain a trigon of type B. 
Hence the reductivity of $P$ is four or less. 
\begin{figure}[ht]
 \begin{center}
  \includegraphics[width=60mm]{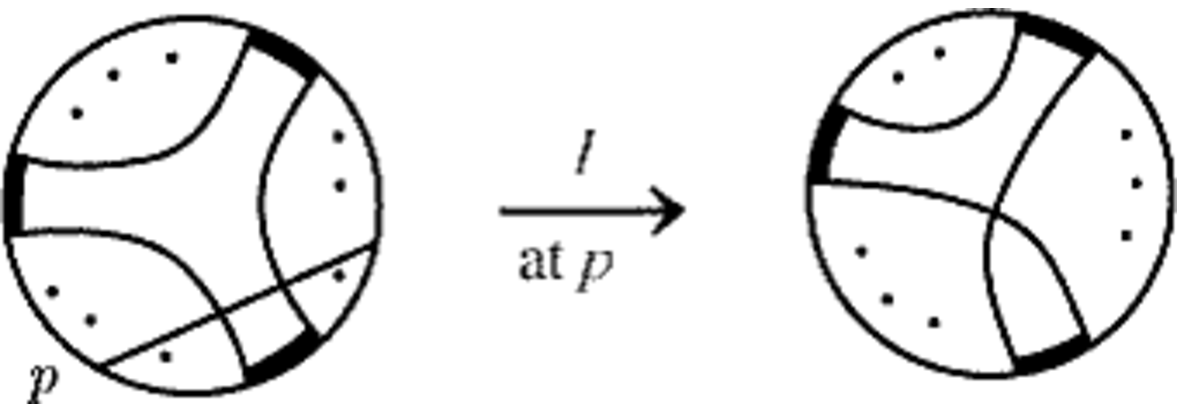}
 \end{center}
 \caption{Trigon of type D.}
 \label{c-3d}
\end{figure}
\end{proof}
\label{c-3d}
\end{lemma}
\phantom{x}

\noindent Adams, Shinjo and Tanaka showed the following lemma \cite{AST}: 

\phantom{x}
\begin{lemma}[Adams-Shinjo-Tanaka]
Every reduced spherical curve has a bigon or trigon. 
\label{ast-lem}
\end{lemma}
\phantom{x}

\noindent They proved this lemma by considering the equality $2C_2 +C_3 =8+C_5 +2C_6 +3C_7+4C_8 +\dots $ which is obtained by the Euler characteristic of $S^2$, where $C_n$ is the number of $n$-gons. 
We prove Theorem \ref{main}. 

\phantom{x}
\noindent {\it Proof of Theorem \ref{main}.} \ 
If $P$ is reducible, then $r(P)=0$. 
If $P$ is reduced, then $r(P)\le 4$ because $P$ has a bigon or trigon. 
\hfill$\square$\\
\phantom{x}

\section{Unavoidable sets for spherical curves}

\begin{figure}[ht]
 \begin{center}
  \includegraphics[width=60mm]{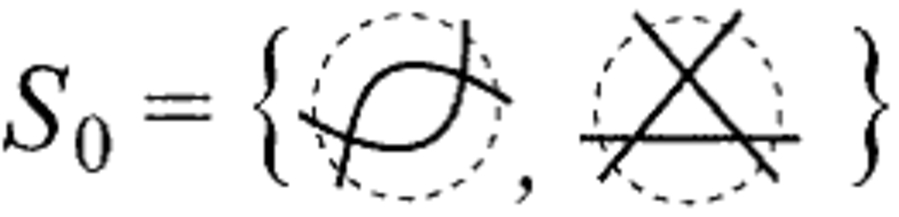}
 \end{center}
 \caption{$S_0$ is an unavoidable set.}
 \label{s1}
\end{figure}

In this section, we discuss the unavoidable sets of tangles for reduced spherical curves. 
A {\it tangle} of a spherical curve is a part of the spherical curve. 
In this section we do not distinguish the mirror image of tangles. 
A set $S$ of tangles is an {\it unavoidable set for spherical curves} if any spherical curve has at least one tangle in $S$. 
For example, the set $S_0$ in Fig. \ref{s1} is an unavoidable set for {\it reduced} spherical curves because of Lemma \ref{ast-lem}. 
We have the following proposition: 

\phantom{x}
\begin{proposition}
The sets $S_1$, $S_2$ and $S_3$ in Fig. \ref{s11} are unavoidable sets for reduced spherical curves. 
\begin{figure}[ht]
 \begin{center}
  \includegraphics[width=130mm]{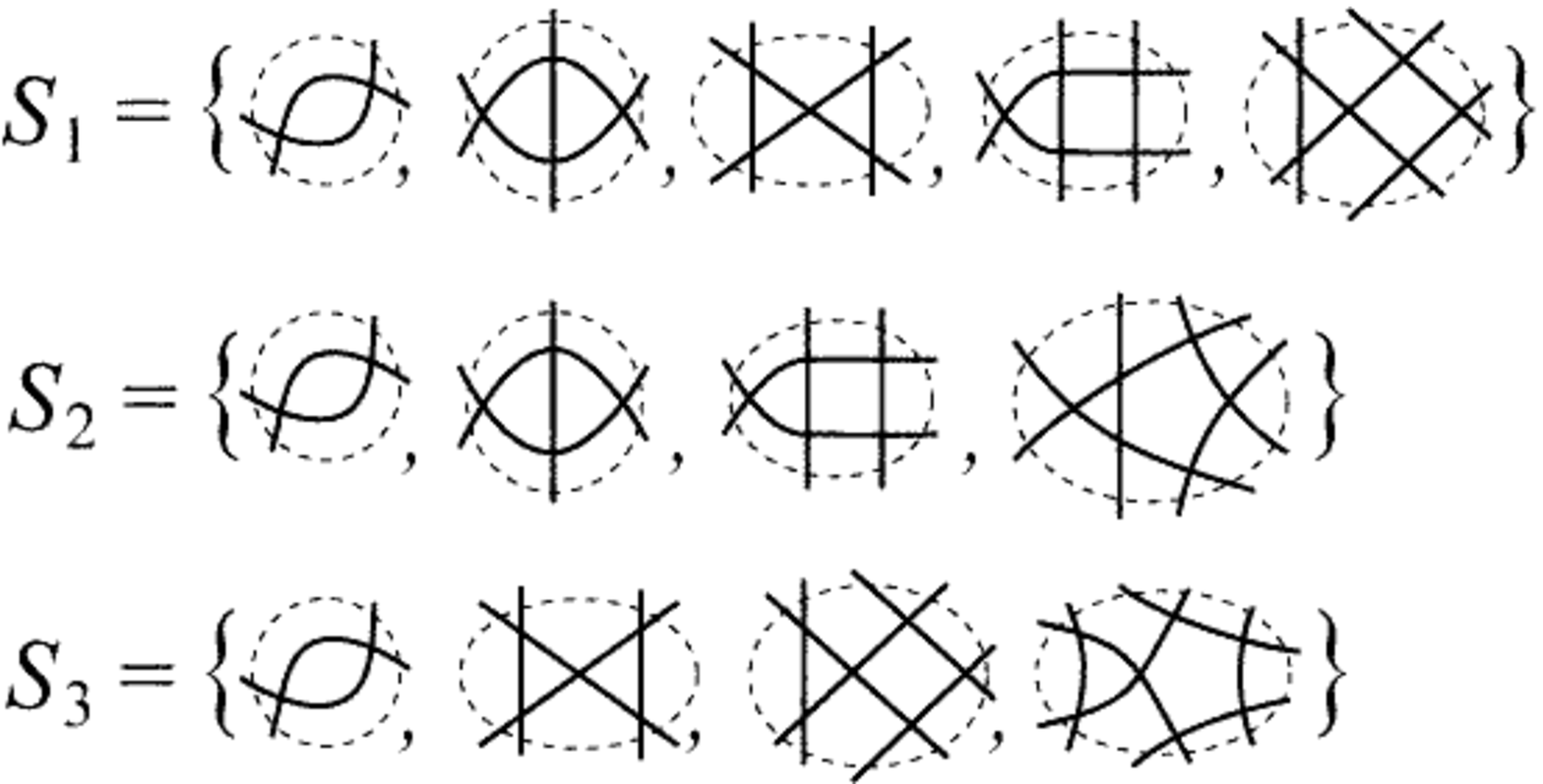}
 \end{center}
 \caption{$S_1$, $S_2$ and $S_3$ are unavoidable sets.}
 \label{s11}
\end{figure}

\begin{proof}
We use the discharging method which is used in graph coloring theory (see, for example, \cite{wilson}). 
We assume there exists a reduced spherical curve $P$ which has none of $S_1$. 
Then $P$ has trigons by Lemma \ref{ast-lem}, and all the regions around each trigon are 5-gons or more. 
We assign a charge of $4-k$ to each $k$-gon. 
Then each trigon receives 1, 4-gon receives 0, 5-gon receives -1, $\dots$. 
The total charge $Ch$ on $P$ is 
$$Ch=C_3-C_5-2C_6-3C_7-4C_8-\dots $$
where $C_n$ is the number of $n$-gons. 
By the equality 
$$2C_2+C_3=8+C_5+2C_6+\dots $$
for Lemma \ref{ast-lem} and that $C_2=0$, we have $Ch=8$. 
Now we do discharging; 
At each trigon, we take off the charge 1, and distribute it to the 6 regions around the trigon by $1/6$ (see Fig. \ref{discharging}). 
\begin{figure}[ht]
 \begin{center}
  \includegraphics[width=70mm]{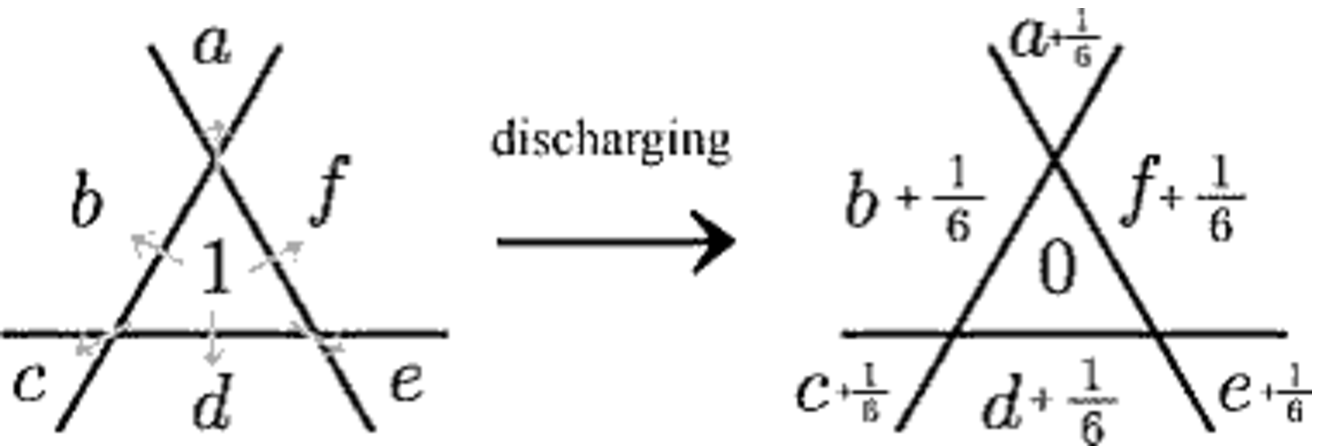}
 \end{center}
 \caption{Discharging.}
 \label{discharging}
\end{figure}
After the discharging, each trigon has the charge 0 because it has been taken the initial charge, and receive no charge from other trigons because they are not adjacent. 
Each 4-gon has also 0 because the initial charge was 0 and receive no charge from trigons. 
Next we consider 5-gons. 
If a 5-gon $R_5$ has $m$ trigons around it, then $R_5$ has the charge $-1+m/6$ after the discharging. 
Remark $m \le 5$ because if $m \ge 6$, some of the trigons around $R_5$ are adjacent. 
Hence each 5-gon has a negative charge after the discharging. 
For $n$-gons ($n=6, 7, 8, \dots$), they have $(4-n)+m/6$ after the discharging, where $m$ is the number of trigons around them. 
We have $(4-n)+m/6<0$ because if $(4-n)+m/6 \ge 0$, then $m \ge 6n-24 \ge 2n$ for $n\ge 6$, and this means all the $2n$ regions around the $n$-gon are trigons, and they are adjacent. 
Therefore, after discharging, the total charge turns a negative number despite that the discharging preserves the total charge 8. 
Hence $S_1$ is an unavoidable set. 
It is similarly shown for $S_2$ and $S_3$ by considering the discharging $D_2$ and $D_3$ in Fig. \ref{discharging2}, respectively. 
\begin{figure}[ht]
 \begin{center}
  \includegraphics[width=130mm]{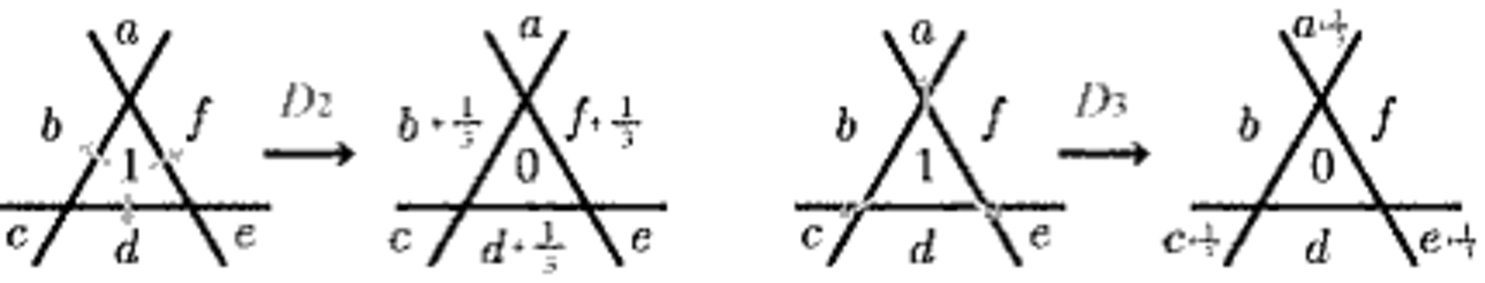}
 \end{center}
 \caption{Discharging $D_2$ and $D_3$.}
 \label{discharging2}
\end{figure}
\end{proof}
\end{proposition}
\phantom{x}

\noindent We have the following question:

\phantom{x}
\begin{question}
Is the set $S_4$ in Fig. \ref{s2} an unavoidable set for reduced spherical curves?
\begin{figure}[ht]
 \begin{center}
  \includegraphics[width=110mm]{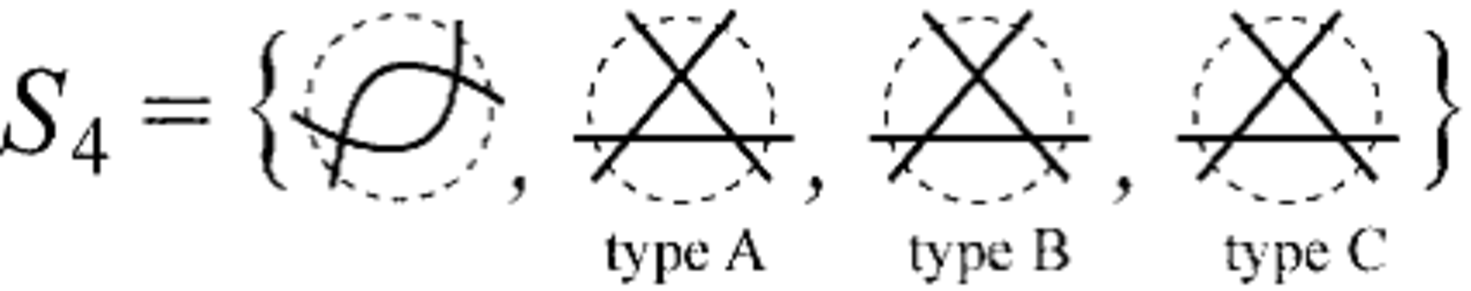}
 \end{center}
 \caption{Unavoidable set?}
 \label{s2}
\end{figure}
\end{question}
\phantom{x}

\begin{figure}[ht]
 \begin{center}
  \includegraphics[width=78mm]{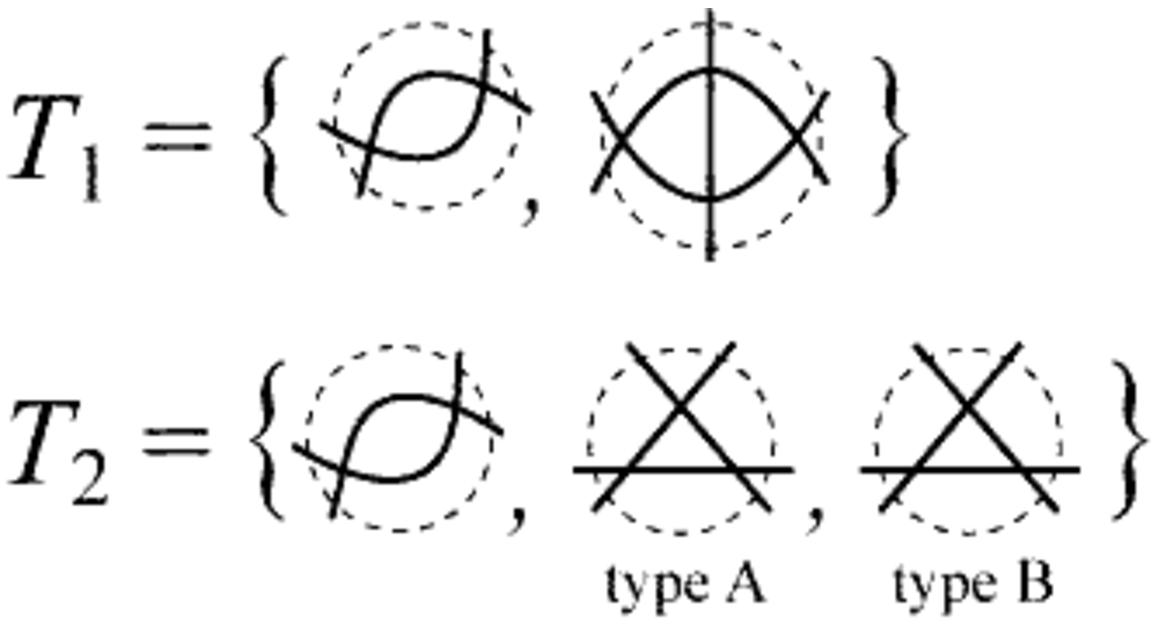}
 \end{center}
 \caption{Avoidable sets.}
 \label{t1}
\end{figure}

\begin{figure}[ht]
 \begin{center}
  \includegraphics[width=90mm]{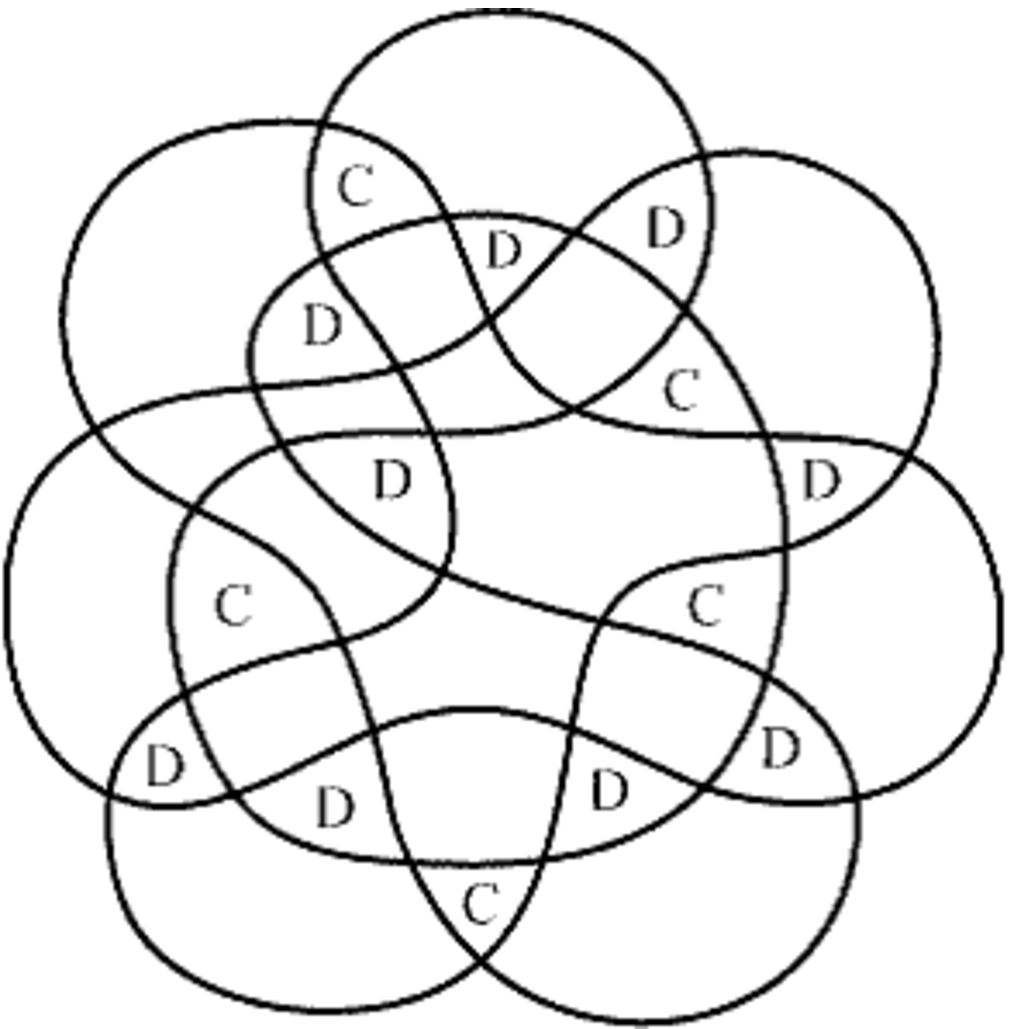}
 \end{center}
 \caption{Spherical curve without bigons or trigons of type A or B.}
 \label{counter}
\end{figure}

\noindent If the answer to this question is Yes, then the answer to Question \ref{main-q} is also Yes. 
A set $T$ of tangles is an {\it avoidable set} for spherical curves if $T$ is not an unavoidable set, i.e., there is a spherical curve which has no tangles in $T$. 
Since we have the spherical curve in Fig. \ref{counter}, the sets $T_1$ and $T_2$ in Fig. \ref{t1} are avoidable sets for reduced spherical curves. 
Remark that $T_2$ is the set consisting of the tangles of the bigon and  the {\it incoherent} trigons. 

\section*{Acknowledgment}

\noindent The author is deeply grateful to Professors Tetsuya Abe, Reiko Shinjo and Kokoro Tanaka for valuable discussions and suggestions. 
She is also grateful to the members of the Friday Seminar on Knot Theory at Osaka City University and Kobe Topology Seminar for helpful advice and encouragements. 
She is supported by Grant-in-Aid for Research Activity Start-up (24840030).

\end{document}